\documentclass[12pt,fleqn]{article}
\usepackage{amsmath,amsfonts,amssymb,amsthm}

\begin{document}
\author{\textbf{Denis~I.~Saveliev}
\thanks{The author was partially supported by RFFI~06-01-00608-a.}}
\title{\textbf{A~Report on a~Game \\on the~Universe of~Sets}}
\date{}
\maketitle

\newtheorem{thm}{Theorem}
\newtheorem{lm}{Lemma}

\newcommand{\cf}{ {\mathop{\mathrm {cf\,}}\nolimits} }
\newcommand{\tc}{ {\mathop{\mathrm {tc\,}}\nolimits} }
\newcommand{\lra}{ {\leftrightarrow} }

\newcommand{\Z}{ {\mathrm Z} }
\newcommand{\ZF}{ {\mathrm {ZF}} }
\newcommand{\ZFA}{ {\mathrm {ZFA}} }
\newcommand{\ZFC}{ {\mathrm {ZFC}} }
\newcommand{\E}{ {\mathrm {AE}} }
\newcommand{\AR}{ {\mathrm {AR}} }
\newcommand{\WR}{ {\mathrm {WR}} }
\newcommand{\AF}{ {\mathrm {AF}} }
\newcommand{\AC}{ {\mathrm {AC}} }
\newcommand{\DC}{ {\mathrm {DC}} }
\newcommand{\AD}{ {\mathrm {AD}} }
\newcommand{\AInf}{ {\mathrm {AInf}} }
\newcommand{\AU}{ {\mathrm {AU}} }
\newcommand{\AP}{ {\mathrm {AP}} }
\newcommand{\APr}{ {\mathrm {APr}} }
\newcommand{\ASp}{ {\mathrm {ASp}} }
\newcommand{\ARp}{ {\mathrm {ARp}} }
\newcommand{\AFA}{ {\mathrm {AFA}} }
\newcommand{\BAFA}{ {\mathrm {BAFA}} }
\newcommand{\FAFA}{ {\mathrm {FAFA}} }
\newcommand{\SAFA}{ {\mathrm {SAFA}} }
\newcommand{\nonempty}{ {\mathrm {nonempty}} }

\newcommand{\I}{ {\mathrm {I}} }
\newcommand{\II}{ {\mathrm {II}} }

\begin{abstract}
In set theory without the~axiom of~regularity,
we consider a~game
in which
two~players choose in~turn
an~element of a~given set,
an~element of this~element, etc.;
a~player wins
if its adversary cannot make any next move.
Sets that are~winning,
i.e.~have a~winning strategy for a~player, 
form a~natural hierarchy with~levels indexed by~ordinals.
We show that
the~class of~hereditarily winning sets
is an~inner model containing all~well-founded sets,
and that
all four~possible relationships between
the~universe,
the~class of~hereditarily winning sets,
and the~class of~well-founded sets
are consistent.
We describe classes of~ordinals
for which it~is consistent that 
winning sets without~minimal elements 
are exactly in the~levels 
indexed by~ordinals of~this class.
For consistency results,
we propose a~new method for getting
non-well-founded models. 
Finally,
we establish a~probability result
by~showing that
on~hereditarily finite well-founded sets
the~first player wins almost always.
\end{abstract}

\par
We work in~ZF, the Zermelo~--- Fraenkel set theory,
minus~AR, the~Axiom of Regularity;
we use~AC, the~Axiom of Choice,
nowhere, for one exception.
Let
$\ZF^-$ and~$\ZFC^-$ be $\ZF$ and~$\ZFC$ minus~AR.
The~notations are standard in~set theory (see~[Jech]):
$V$~is the~universe,
i.e.~the~class of sets,
$On$~is the~class of~ordinals,
$V_\alpha$~is the~$\alpha$th level of the~cumulative hierarchy,
$\tc$\!~is the~transitive closure,
$\cf$\!~is the~cofinality,
etc.
Let also
$V_\infty$~be the~class of well-founded sets
and $r$~the~rank function.

\par
Suppose 
two players, I and~II,
starting from a~given set~$x$,
try to~construct
an~$\in$-decreasing sequence of sets~$x_n$
such that
the~adversary cannot continue~it:
first I chooses $x_0\in x$,
then II chooses $x_1\in x_0$, etc.;
a~player {\it wins\/} 
if its~adversary cannot make any next move,
i.e.~if he could choose the~empty set.
(This game~--- but not our results~--- can be found in~[For]
where it~is considered in~Quine's NF
(the~accont has little mistakes);
a~close game is in~[BaMo].)
We say that a~set is~{\it winning\/}
if it has a~winning strategy for a~player.
Let
$W$~be the~class of~winning sets,
$W_\I$ and~$W_\II$
the~classes of sets winning for~I and~II~resp.
Clearly,
$W_\I\cup W_\II=W$
and
$W_\I\cap W_\II=0$.
The~class $W$~admits 
a~natural hierarchy:
Given an~ordinal~$\gamma$,
a~set is
$2\gamma$-{\it winning\/}
if each of~its elements is $2\delta+1$-winning
for some~$\delta<\gamma$,
and
$2\gamma+1$-{\it winning\/}
if some its element is~$2\gamma$-winning.
Let
$W_\nu$ be the~class of $\nu$-winning sets.
Clearly,
$W_0=\{ 0\}$,
$
W_{2\gamma+1}=
V\setminus
{\cal P}(V\setminus W_{2\gamma})
$,
$
W_{2\gamma}=
{\cal P}(\bigcup_{\delta<\gamma}W_{2\delta+1})
$,
and
$W=\bigcup_{\nu\in On}W_\nu$.
Also
$
W_\I=
V\setminus{\cal P}(V\setminus W_{\II})=
\bigcup_{\gamma\in On} W_{2\gamma+1}
$
and
$
W_\II=
{\cal P}(W_\I)=
\bigcup_{\gamma\in On}W_{2\gamma}
$.
Let
$S_\nu$ be $\nu$th level in the~winning hierarchy:
$S_\nu=W_\nu\setminus\bigcup_{\mu<\nu}W_\mu$.
Clearly,
$W_{2\gamma}=
\bigcup_{\delta\le\gamma}S_{2\delta}$
and
$W_{2\gamma+1}=
\bigcup_{\delta\le\gamma}S_{2\delta +1}$.
The~classes $S_\nu$ can be defined recurrently:
$$
\begin{array}{lll}
S_{2\gamma+1}=
{\cal P}
(V\setminus
\bigcup_{\delta<\gamma}S_{2\delta})
\,\setminus\,
{\cal P}
(V\setminus
S_{2\gamma}),\\
S_{2\gamma}=
{\cal P}
(\bigcup_{\delta<\gamma}S_{2\delta+1})
\,\setminus\,
\bigcup_{\varepsilon<\gamma}
{\cal P}
(\bigcup_{\delta<\varepsilon}S_{2\delta+1}).
\end{array}
$$
For $x\in W$
we let
$w(x)=\min\{\nu:x\in S_\nu\}$.
If~$x$~is nonempty,
$$
w(x)=
\left\{
\begin{array}{lll}
\min                          
\{w(x_0)+1:x_0\in x\cap W_\II\}
&\text{if }\;
x\in W_\I,\\
\sup
\{ w(x_0)+1:x_0\in x\}
&\text{if }\;
x\in W_\II.
\end{array}
\right.
$$

\begin{lm}
$S_\nu$ is nonempty
for all~$\nu$.
\end{lm}

Let
$HW$ consist 
of hereditarily winning sets:
$HW=\{x:\tc(\{x\})\subseteq W\}$.

\begin{lm}
${\cal P}(W)\subseteq W$
and\,
${\cal P}(HW)=HW$.
\end{lm}

\begin{thm}
$HW$ is an~inner model (of~$\ZF^-$).
Moreover,
$HW\supseteq V_\infty$.
\end{thm}

\par
Which from the~classes
$V\supseteq W\supseteq HW\supseteq V_\infty$
can be distinct?
Clause~(1) of Lemma~3 says that
it suffices to~know this for 
$V\supseteq HW$ and $HW\supseteq V_\infty$;
Theorem~4 below shows that
all four remaining cases are consistent.

\begin{lm}
\rule{0mm}{0mm}
\\
1.
$V=W$ is equivalent to $W=HW$.
\\
2.
If\,
$V\setminus W$ is nonempty
\,then
$V\setminus W$ 
\!and\,
$W\setminus HW$\! 
are proper.
\\
3.
If\,
$HW\setminus V_\infty$ is nonempty
\,then
it is proper.
\\
4.
$\AR$~is equivalent to~$W=V_\infty$.
\end{lm}

\par
Given~$\nu$,
how many $\nu$-winning well-founded sets 
have a~given rank?
Let
$S_{\alpha,\nu}=V_\alpha\cap S_\nu$,
or recurrently,
$$
\begin{array}{lll}
S_{\alpha+1,2\gamma+1}=
{\cal P}
(V_\alpha\setminus
\bigcup_{\delta<\gamma}S_{\alpha,2\delta})
\,\setminus\,
{\cal P}
(V_\alpha\setminus
S_{\alpha,2\gamma}),\\
S_{\alpha+1,2\gamma}=
{\cal P}
(\bigcup_{\delta<\gamma}S_{\alpha,2\delta+1})
\,\setminus\,
\bigcup_{\varepsilon<\gamma}
{\cal P}
(\bigcup_{\delta<\varepsilon}S_{\alpha,2\delta+1}).
\end{array}
$$

\begin{lm}
\rule{0mm}{0mm}
\\
1.
$w(x)\le r(x)$
for all well-founded $x$.
\\
2.
$S_{\alpha,\nu}$ is nonempty
\,iff\;
$\nu<\alpha$.
\\
3.
$
S_{\alpha+1,\nu}
\setminus
S_{\alpha,\nu}
$
is nonempty
\,iff\;
$0=\nu=\alpha$
or
$0<\nu\le\alpha$.
\\
4.
$
S_{\alpha,\nu}\in
S_{\alpha+1,\nu+1}
\setminus
S_{\alpha,\nu+1}
$
\,iff\;
$0=\nu<\alpha=1$
or
$0<\nu<\alpha$.
\end{lm}

\begin{lm}
In~$s_{\alpha,\nu}$ below,
let
$\alpha>\nu$
and, except case~(1),~$\nu>0$.
Then:
\\
1. (Trivial case)
For all $\alpha$,
$
|S_{\alpha,0}|=1
$.
\\
2. (Finite case)
For all $m<\omega$,
$$
\begin{array}{lll}
|S_{m+1,2k+1}|=
2^{|V_m|-\sum_{j<k}|S_{m,2j}|}
-
2^{|V_m|-\sum_{j\le k}|S_{m,2j}|},\\
|S_{m+1,2k+2}|=
2^{\sum_{j\le k}|S_{m,2j+1}|}
-
2^{\sum_{j<k}|S_{m,2j+1}|}.
\end{array}
$$         
\\
3. (Infinite case)
For all
$\alpha\ge\omega$,
$
|S_{\alpha,\nu}|=
|V_{\alpha}|
$
and\,
$
|S_{\alpha+1,\nu}\setminus S_{\alpha,\nu}|=
|V_{\alpha+1}|
$.
\end{lm}

\par
This leads to an~interesting probability result:
Let
$$
\Pr(S_{\omega,n})=
\lim_{m\to\omega}
\frac{|S_{\omega,n}\cap V_m|}{|V_m|}.
$$
Easy calculations show that
the~limits are exist
(and $\Pr$ generates a~probability on~$V_\omega$);
moreover,

\begin{thm}
$
\Pr(S_{\omega,1})=
\Pr(S_{\omega,3})=
1/2
$
\,(and
\,$\Pr(S_{\omega,n})=0$
for all other~$n$).
\end{thm}

\par
Thus
a~``half'' of hereditarily finite well-founded sets
consists of 1-winning sets,
while another ``half'' consists of 3-~but not~1-winning ones.
It follows that
almost always player~I has a~winning strategy;
moreover, he may win very quickly:
at~1 or~3 moves.

\par
How similar is~$W$ to a~model of~$\ZFC$?
It can be shown that
either $W=V$
or $W$ cannot hold many ZFC~axioms.
Surprisingly,
Regularity is not from them:
$\AR^W$\!~is consistent with~$\neg\AR$.
Let $\varsigma(C)$ means
``some $x\in C$ has no~$\in$-minimal elements''.

\begin{lm}
\rule{0mm}{0mm}
\\
1.
$\neg\AR^W$\! is equivalent to
\,$\varsigma({\cal P}(W))$.
\\
2.
$\neg\AR^{HW}$\! is equivalent to
\,$\varsigma(HW)$.
\\
3.
$\AR^{HW}$\! is equivalent to
$HW=V_\infty$.
\\
4.
$\AR^W$\! implies $\AR^{HW}$\!.
\end{lm}

\par
For consistency results,
we propose a~new construction of non-well-founded models
(see also~[Sav]).
A~customary procedure
(based on~quotients under bisimulations;
see [FoHo], [Acz], [BaMo], [d'Ag])
leads to~models too ``unstratified'' for~us.
Our~construction is rather like a~cumulative hierarchy:
we put an~appropriate~$(M_0,E_0)$ as an~initial level
and apply iteration of a~certain analog of power set operation.
In the~derived hierarchy,
upper levels end-extend lower ones;
interesting properties are reflected at~lower levels
and so easily controlled.
We say that
$v\in A$
{\it represents\/} a~set $y\subseteq A$ in~$(A,R)$
if
$y=\{ u:u\,R\,v\}$,
and that
$(A,R)$ is {\it thick\/}
if
for every set $y\subseteq A$
there exists $v\in A$
such that
$v$ represents~$y$ in~$(A,R)$.
Put
${\cal S}(A)\subseteq{\cal P}(A)$
consist of all nonempty subsets of~$A$
except for the~represented in~$(A,R)$.
Let
$B\supseteq A$.
$(B,S)$ {\it end-extends\/}~$(A,R)$
if
$S\cap(B\times A)=R$.
$(A,R)$ {\it reflects\/} a~formula~$\varphi $ over~$(B,S)$
if
$\varphi^{A,R}(x,\ldots)\,\lra\,\varphi^{B,S}(x,\ldots)$
for all~$x,\ldots\in A$.
Given
$(M_0,E_0)$,
we define
$(M_\alpha,E_\alpha)$
by recursion on~$\alpha$:
$$
\begin{array}{rll}
M_{\alpha+1}=
M_\alpha\cup {\cal S}(M_\alpha),
&E_{\alpha+1}=
E_\alpha\cup\,(\in\cap\,(M_\alpha\times{\cal S}(M_\alpha))),\\
M_\alpha=
\bigcup_{\beta<\alpha}M_\beta,
&E_\alpha=
\bigcup_{\beta<\alpha}E_\beta
\;\;\;\text{ if $\alpha$ is limit},\\
M=
\bigcup_{\alpha\in On}M_\alpha,
&E=
\bigcup_{\alpha\in On}E_\alpha.
\end{array}
$$

Really,
we need $(M_0,E_0)$
satisfying the~properties~(i)-(iv):
$$
\begin{array}{rll}
\text{ (i) }
&\{ u:u\,E_0\,v\},\\
\text{ (ii) }
&x\notin\tc(y)
\;\text{ if }x,y\in M_0,\\
\text{ (iii) }
&(M_0,E_0)\models\text{ Empty~Set},\\
\text{ (iv) }
&(M_0,E_0)\models\text{ Extensionality}.
\end{array}
$$

\begin{lm}
Suppose
$(M_0,E_0)$ satisfies {\em (i)}-{\em (ii)}.
Then:
\\
1.
$(M,E)$ end-extends all the~$(M_\alpha,E_\alpha)$.
\\
2.
$(M_0,E_0)$
reflects
Extensionality over~$(M,E)$.
\\
3.
If
$(M_0,E_0)$ satisfies also~{\em (iii)},
then
$(M,E)$~is thick.
\\
4.
$(M_1,E_1)$
reflects
\,$\varsigma(W_{2\gamma})$
and
\,$\varsigma(W_\II)$
over~$(M,E)$.
\\
5.
$(M_1,E_1)$
reflects\,
$\AR^W$\!
over~$(M,E)$.
\\
6.
$(M,E)\models\varsigma(S_{2\gamma})$
implies
$(M_1,E_1)\models\varsigma(S_{2\delta})$
for some~$\delta\le\gamma $.
\end{lm}

We need~(4)-(6) only for results on winning sets;
by a~theorem in~[Rieg],
(1)-(3) suffice to~hold $\ZF^-$~axioms:

\begin{thm}
If
$(M_0,E_0)$
satisfies~{\em (i)-(iv)},
then
$(M,E)$ is a~model of~$\ZF^-$.
Moreover, 
it is a~model of~$\ZFC^-$ (if so is~$V$).
\end{thm}

\begin{thm}
All the consistent relationships between
$V$, $W$, $HW$, $V_\infty$,
and relativizations of~$\AR$ to~them
are exactly these five:
\\
1.
$
V=W=HW=V_\infty
+\AR+\AR^W+\AR^{HW}
$,
\\
2.
$
V\ne W\ne HW=V_\infty
+\neg\AR+\AR^W+\AR^{HW}
$,
\\
3.
$
V\ne W\ne HW=V_\infty
+\neg\AR+\neg\AR^W+\AR^{HW}
$,
\\
4.
$
V=W=HW\ne V_\infty
+\neg\AR+\neg\AR^W+\neg\AR^{HW}
$,
\\
5.
$
V\ne W\ne HW\ne V_\infty
+\neg\AR+\neg\AR^W+\neg\AR^{HW}
$.
\end{thm}

Next we have
more fine results on~$\varsigma(C)$
for various~$C\subseteq W$:

\begin{lm}
\rule{0mm}{0mm}
\\
1.
$\varsigma (W_\I)$
is equivalent to\,
$\neg\AR$.
\\
2.
$\varsigma (W_\II)$
implies\,
$\neg\AR^W$\!.
\\
3.
$\neg\varsigma(W_\II)$
is consistent with\,
$\neg\AR^W$\!.
\\
4.
$\varsigma(W_\II)$
is consistent.
\\
5.
$\neg\varsigma({\cal P}(W_\II))$.
\end{lm}

\begin{lm}
\rule{0mm}{0mm}
\\
1.
$\varsigma(S_\nu)$
\,implies\,
$\nu>1$.
\\
2.
$\neg\AR$
\,implies\,
$\varsigma(S_\nu)$
for all odd~$\nu>1$.
\\
3.
If 
$\mu$ is even,
$\,\varsigma(S_\mu)$
\,implies\,
$\varsigma(S_\nu)$
for all~$\nu\ge\mu$.
\\
4.
If 
$\mu$ is limit,
$\,\varsigma(S_\mu)$
\,implies\,
$\varsigma(S_\nu)$
for some~$\nu\le\mu$ with~$\cf\nu=\omega$.
\end{lm}

\begin{thm}
Let
$A\subseteq On$.
Then 
$\{\nu:\varsigma(S_\nu)\}=A$~is consistent 
with~$\ZFC^-$
exactly in one from three cases:
\\
1.
$A$ is empty,
\\
2.
$A=\{\nu>1:\nu\text{ is odd\/}\}$,
\\
3.
$A=\{\nu>1:\nu\text{ is odd or }\nu\ge\mu\}$
for some $\mu$ with $\cf\mu\le\omega$.
\end{thm}

\par
The~only use of~AC is the~proof of~(4)~of Lemma~9
(and then (3)~of Theorem~5);
I~do not know can it be omitted.
The~consistency results remain valid
under assumptions which add no~new winning sets, as
``all sets are strongly extensional'' 
(see [Acz],~[BaMo]),
``there are no~$\in$-cyc\-les'',
``there are hereditarily non-winning sets of every structure'',
``every set without $\in$-mini\-mal elements has size greater than any given'',
etc.
Summary,
we see a~deep difference
between odd-~and even-winning sets
by finding the~latters more rare and queer.

\end{document}